\input amstex
\documentstyle{amams} 
\document
\annalsline{157}{2003}
\received{May 16, 2001}
\startingpage{545}
\font\tenrm=cmr10
\input amssym.def
\input amssym.tex

\catcode`\@=11
\font\twelvemsb=msbm10 scaled 1100

\font\ninemsb=msbm10 scaled 800
\newfam\msbfam
\textfont\msbfam=\twelvemsb  \scriptfont\msbfam=\ninemsb
  \scriptscriptfont\msbfam=\ninemsb
\def\msb@{\hexnumber@\msbfam}
\def\Bbb{\relax\ifmmode\let\next\Bbb@\else
 \def\next{\errmessage{Use \string\Bbb\space only in math
mode}}\fi\next}
\def\Bbb@#1{{\Bbb@@{#1}}}
\def\Bbb@@#1{\fam\msbfam#1}
\catcode`\@=12

 \catcode`\@=11
\font\twelveeuf=eufm10 scaled 1100
\font\teneuf=eufm10
\font\nineeuf=eufm7 scaled 1100
\newfam\euffam
\textfont\euffam=\twelveeuf  \scriptfont\euffam=\teneuf
  \scriptscriptfont\euffam=\nineeuf
\def\euf@{\hexnumber@\euffam}
\def\frak{\relax\ifmmode\let\next\frak@\else
 \def\next{\errmessage{Use \string\frak\space only in math
mode}}\fi\next}
\def\frak@#1{{\frak@@{#1}}}
\def\frak@@#1{\fam\euffam#1}
\catcode`\@=12

\title{On a coloring conjecture\\ about unit fractions}
\shorttitle{On a coloring conjecture  about unit fractions} 
  \author{Ernest S.\ Croot III}
   \institutions{University of California, Berkeley, CA\\
 {\eightpoint {\it E-mail address\/}: ecroot@math.berkeley.edu}}

 \vfil
 \centerline{\bf Abstract}
\vglue12pt
We prove an old conjecture of Erd{\H o}s and Graham on sums of unit 
fractions:  There exists a constant $b>0$ such that if we $r$-color the 
integers in
$[2,b^r]$, then there exists a monochromatic set $S$ such that 
$\sum_{n \in S} 1/n = 1$.
\vglue18pt
\section{Introduction}

We will prove a result on unit fractions which has the following corollary. 

\nonumproclaim{Corollary} 
There exists a constant $b$ so that for every partition of the integers in
$[2, b^r]$ into $r$ classes{\rm ,}
there is always one class containing a subset $S$ with the property
$\sum_{n \in S} 1/n = 1$.  
\endproclaim

In fact, we will show that $b$ may be taken to be $e^{167000}$, if 
$r$ is sufficiently large, though we believe that $b$ may be taken to be 
much smaller;  also  note that $b$ cannot be taken to be smaller than $e$, 
since the integers in $[2,e^{r - o(r)}]$ can be placed into $r$ classes in 
such a way that the sum of reciprocals in each class is just under $1$.

This corollary implies the result mentioned in the abstract and so resolves an
unsolved problem of Erd{\H o}s and Graham, which appears in [2], [3], and
[5].

We will need to introduce some notation and definitions in order to state the
Main Theorem, as well as the propositions and lemmas in later sections:  
For a given set of integers $C$, let $\Cal{Q}_C$ denote the set 
of all the prime power divisors of elements of $C$, and let 
$\Sigma(C) = \sum_{q \in \Cal{Q}_C} 1/q$.
Define $\Cal{C}(X, Y; \theta)$ to be the integers in $[X,Y]$ all of whose 
prime power divisors are $\leq X^\theta$, and
let $\Cal{C}'(X, Y; \theta)$ be those integers  
$n \in \Cal{C}(X, Y; \theta)$ such that $\omega(n) \sim \Omega(n) \sim 
\log\log n$, where $\omega(n)$ and $\Omega(n)$
denote the number of prime divisors and the number of prime power divisors
of $n$, respectively. 

Our Main Theorem, then, is as follows.
\nonumproclaim{Main Theorem}   Suppose $C \subset \Cal{C}'(N, N^{1+ \delta};
\theta)${\rm ,} where $\theta, \delta > 0$, and $\delta + \theta < 1/4$.  
If $N \gg_{\theta,\delta} 1$ and  
$$
\sum_{n\in C} {1 \over n} >  6,
$$
then there exists a subset $S\subset C$ for which $\sum_{n\in S} 
1/n =1$. 
\endproclaim

To prove the corollary, we will show in the next section that for $r$
sufficiently large,
$$
\sum_{n \in \Cal{C}'(N,N^{1+\delta}; 1/4.32)} {1 \over n}\ >\ 6r, \tag{1.1}
$$
where $N = e^{163550r}$ and $N^{1+\delta} = e^{166562r}$.  Thus, if we 
partition the integers in $[2,e^{167000r}]$ into $r$ classes,
then for $r$ sufficiently large, one of the classes $C$ satisfies the 
hypotheses of the Main Theorem, and so our corollary follows.

The key idea in the proof of the Main Theorem is to construct a subset of 
$C$ with usable properties. These are summarized in the following 
proposition which is proved in Section 4.  

\nonumproclaim{Proposition 1}  Suppose 
$C \subset \Cal{C}'(N, N^{1+\delta}; \theta)$ with
$\delta + \theta < 1/4${\rm ,} and suppose
$$
\sum_{n \in C} {1 \over n}\ >\ 6.
$$
Then there exists a subset $D \subset
C$ such that  
$$
\sum_{n \in D} {1 \over n} \in [2 - 3/N, 2),\tag{1.2}
$$
and which has the following property\/{\rm :}  If $I$ is an interval of length
$N^{3/4}$ for which there are less than $N^{1-\theta}/(\log\log N)^2$ 
elements of $D$
that do not divide any element of $I${\rm ,} then every element of $D$
divides one single element of $I$.
\endproclaim

The sum of the reciprocals of the elements of $D$ is $<2$ by (1.2),
so if there is a subset $S$ of $D$ for which 
$\sum_{n\in S} 1/n$ is an integer then that sum equals 1 or $S$ is the 
empty set. Now if $x$ is an integer and
$$
P\ :=\ \text{\rm lcm} \{ n \in D\},
$$
then $(1/P)\sum_{h\pmod P} e( hx/P)=1$ if $x/P$ is an integer, 
and is $0$ otherwise, where $e(t)=e^{2\pi i t}$. 
Combining these remarks we deduce \pagebreak that  
$$
\#\left\{ S \subset D:\ \sum_{n\in S} 1/n =1 \right\}=\ \left ( {1 \over P} \sum_{-P/2 < h \leq P/2} E(h) \right ) - 1, \tag{1.3}
$$
where 
$$
E(h)\ :=\ \prod_{n \in D} \left ( 1 + e(h/n) \right ).
$$

Now,
$$
E(h)\ =\ e\left ( {h \over 2} \left \{ \sum_{n \in D} {1 \over n}
\right \} \right ) \left ( 2^{|D|} \prod_{n \in D}
\cos ( \pi h/n ) \right ), \tag{1.4}
$$
so that 
$$
\text{Arg}(E(h))\ =\ \pi h \left \{ \sum_{n \in D} {1 \over n}
\right \} \in (2\pi h - \pi/2, 2\pi h + \pi/2),
$$
if $|h|$ is an integer $< N/6$; and therefore
$E(h)+E(-h)>0$ for this case.  Thus we deduce that 
$$
\sum_{|h| < N/6} E(h) > E(0) = 2^{|D|}.
$$

For $h$ in the range $N/6 \leq |h| \leq P/2$, we will use Proposition 1 to
show that
$$
|E(h)| < {2^{|D|-1} \over P},\tag{1.5}
$$
so that, by the last two displayed equations, 
$$
{1 \over P} \sum_{-P/2 < h \leq P/2} E(h) > {1 \over P} \left (
2^{|D|} - \sum_{ N/6 \leq |h| \leq P/2} {2^{|D| - 1} \over 
P} \right )> {2^{|D|-1}  \over P} > 1,
$$
since $|D|\geq \sum_{n\in D} N/n \geq 2N-3$, and since 
$$
P < \left ( N^\theta\right )^{\pi(N^\theta)} \ll e^{(1+o(1))N^\theta} = o\left (
2^{|D|}\right ), \tag{1.6}
$$
by the prime number theorem.  Theorem 1 then follows.  

We will now see how (1.5) follows from Proposition 1.  
If $|h| \in [N/6, P/2]$ then $I := [h-N^{3/4}/2,
h+N^{3/4}/2]$ does not contain any integer divisible by every element of
$D$, since $P = \text{lcm}_{n \in D} n$ is bigger than every element in $I$.  
Therefore, by Proposition 1, there are at least 
$N^{1-\theta-o(1)}$ elements 
$n \in D$ which do not divide any integer in $I$.  For such $n$ 
we will have that $||h/n||> N^{3/4}/2n > 1/(2 N^{1/4+\delta})$ 
(where $||t||$ denotes the distance from $t$ to the nearest integer to $t$).  
Thus, 
$$
\eqalign{
\left | \prod_{n \in D} \cos(\pi h/n) \right |\  &< \  
\left | \cos\left ({\pi \over 2N^{1/4 + \delta}}\right ) 
\right |^{ N^{1-\theta-o(1)}}\cr 
&< \left( 1 - {\pi^2 \over 8N^{1/2 + 2\delta}} + O \left (
{1 \over N}\right ) \right )^{ N^{1-\theta-o(1)} }\cr
&< \exp\left ( - (\pi^2/8)  N^{1/2-2\delta - \theta-o(1)}\right ) 
< {1 \over 2P},
}
$$
by (1.6) since $\delta + \theta < 1/4$, and so (1.5) follows from (1.4).  

The rest of the paper is dedicated to proving Proposition 1.

\section{Normal integers with small prime factors}

We will need the following result of Dickman from [1].
\nonumproclaim{Lemma 1} Fix $u_0 > 0$.  For any $u${\rm ,}  $0 < u < u_0$ we have
$$
\#\{ n \leq x : p|n \Rightarrow p \leq x^{1/u}\} \sim x\rho(u),
$$  
where $\rho(u)$ is the unique{\rm ,} continuous solution to
the differential difference equation
$$
\cases \rho(u) = 1, &\text{if}\  0 \leq u \leq 1\\
       u\rho'(u) = -\rho(u-1), &\text{if}\  u > 1.
\endcases
$$
\endproclaim

From this lemma and partial summation we have, for a fixed $u$ and $\delta$, 
$$
\sum_{N < n < N^{1+\delta} \atop p^a || n \Rightarrow 
p^a \leq N^{1/u}} {1 \over n} \sim 
{\log N \over u} \int_u^ {u(1+\delta)} \rho(w) \text{\rm d}w.
$$ 
Using this, a numerical calculation shows for 
$N = \exp(163550r)$, $\theta = 1/u = 1/4.32$, and $\delta = 1/4 - \theta
- 0.0001$ that 
$$
\sum_{N < n < N^{1+\delta} \atop p^a || n \Rightarrow p^a \leq  N^\theta} 
{1 \over n} > 6.0001r. 
$$
Combining this with the well-known fact that almost all integers $n\leq x$ 
satisfy $\omega(n) \sim \Omega(n) \sim \log\log n$, so that 
$$
\sum_{N < n < N^{1+\delta} \atop \omega(n)\text{ or }\Omega(n) \not 
\sim 
\log\log N} {1 \over n} = o(r),
$$
we have that (1.1) follows. 

\section{Technical lemmas and their proofs}
\nonumproclaim{Lemma 2}  If $w_1$ and $w_2$ are distinct integers which both lie 
in an interval of length $\leq N${\rm ,} then 
$$
\sum_{p^a | \text{\rm gcd}(w_1,w_2)} {1\over p^a} < \sum_{p | 
\text{\rm gcd}(w_1,w_2)} {1 \over p} + O(1) < (1+o(1))\log\log\log N.
$$
\endproclaim

\demo{Proof of Lemma {\rm 2}}  Let $G = $gcd$(w_1,w_2)$.
We have that $G \leq |w_1 - w_2| < N$, since $G | |w_1 - w_2|$;
also, $\omega(G) = o(\log N)$, since $\omega(n) = o(\log N)$ uniformly for 
$n \leq N$.  Now, by the Prime Number Theorem,  
we have $\pi(\log N\log\log N) \gg \log N> w(G)$, for $N$ sufficiently
large, and so
\vglue12pt
\hfill ${\displaystyle
\sum_{p | G \atop p\text{ prime}} {1 \over p} < \sum_{p \leq \log N\log\log N
\atop p\text{ prime}} {1 \over p} < (1 + o(1))\log\log\log N.
}$
\enddemo
  \vglue2pt

\nonumproclaim{Lemma 3}  If $H \subset \Cal{C}(N,N^{1+\beta}; 1)${\rm ,} $\beta > 0${\rm ,}
satisfies $$\sum_{n \in H} 1/n >  1/ (\log N)^{o(1)},$$  and $\omega(n) 
\sim \log\log n$, for every $n \in H${\rm ,} then
$$
\Sigma(H) > (e^{-1} - o(1))\log\log N.
$$
\endproclaim

{\it Proof of Lemma {\rm 3}}.  From the hypotheses of the lemma, together with
the fact that $t! > (t/e)^t$ for $t \geq 1$, we have that 
$$
\eqalign{
{1 \over (\log N)^{o(1)}} 
&< \sum_{n\in H} {1 \over n} < \sum \Sb n \ :\  p^a | n \Rightarrow
p^a \in \Cal{Q}_H \\ \omega(n) \sim \log\log n \sim \log\log N 
\endSb {1 \over n}
< \sum_{t \sim \log\log N} {\Sigma(H)^t \over t!}\cr
&< \sum_{t \sim \log\log N} \left ( {\Sigma(H) e \over t}\right )^t
= \left ( {\Sigma(H) (e+o(1)) \over \log\log N} \right )^{(1+o(1))
\log\log N},
}
$$   
and so $\Sigma(H)$ satisfies the conclusion to Lemma 3.
\hfill\qed

\section{Proof of Proposition 1}

Before we prove Proposition 1, we will need two more propositions. 
 
\nonumproclaim{Proposition 2} \hskip-6pt Suppose that 
$J \subset \Cal{C}(N, \infty; \theta)${\rm ,} where $\theta < 1${\rm ,}
and $\sum_{n \in J} 1/n\break \geq \alpha > \nu$.
If $N \gg_{\alpha, \nu,\theta } 1${\rm ,} then there is a subset $E
\subset J$ such that 
$$
 \sum_{n \in E} {1 \over n} \in \left [ \nu - {1 \over N}, \nu
\right );\tag{4.1}
$$  and{\rm ,}
$$ \sum_{n \in E \atop q|n} {1 \over n} > { \min\{ \nu, \alpha -
\nu\} \over 5q \log\log N},
\text{ for all $q \in \Cal{Q}_E$}. \tag{4.2}
 $$
\endproclaim

\nonumproclaim{Proposition 3}  Suppose that $E\subset \Cal{C}'(N,N^{1+\delta}; \theta)${\rm ,} 
$0 < \theta < 1/4${\rm ,} satisfies {\rm (4.1)} and {\rm (4.2).} 
If all but at most $N^{1-\theta}/(\log\log  N)^2$ 
elements of $E$ divide some element of an interval 
$I := [h - N^{3/4}/2, h + N^{3/4}/2]${\rm ,}  then either\medbreak

\item{\rm A.}  There is a single integer in $I$ divisible by all elements of $E${\rm ,} or
\vglue4pt
\item{\rm B.}  There exist distinct integers $w_1, w_2 \in I${\rm ,} such that
\vglue-12pt
$$
\#\{ n \in E\ :\ n \nmid w_1  \text{ and } n \nmid w_2 \} < 
{2N^{1-\theta} \over (\log\log N)^2 }, \tag{4.3}
$$
$$
\text{\rm lcm} \{ n \in E\} = \text{\rm lcm} \{q \in \Cal{Q}_E\}
 | w_1 w_2, \tag{4.4}$$
\item{}and
\vglue-12pt
$$ \align
(e^{-1} - o(1))\log\log N &< \sum_{q | w_i \atop q\in \Cal{Q}_E} {1 
\over q}\tag{4.5}
\\
& < (1 - e^{-1} + o(1))\log\log N,\text{ for $i=1$ and $2$}.
\endalign $$
\endproclaim  

These propositions will be proved in the next two sections of the paper.
To prove Proposition 1, we iterate the following procedure:\medbreak

1.  Set $j =  0$ and let $C_0 := C$.
\vglue4pt
2.  Use Proposition 2 with $J = C_j$, $\alpha = \sum_{n \in
C_j} 1/n > 2$ , and $\nu = 2$, to produce a subset $E$ satisfying 
(4.1) and (4.2).
\vglue4pt
3.  If case A of Proposition 3 holds for every real number $h$ satisfying the 
hypotheses of Proposition 3, then we can let $D := E$, and
Proposition 1 is proved.
\vglue4pt
4.  If there is some $h$ for which case B holds, then, by (4.3), we have 
for either $i=1$ or $i=2$ that 
$$ \align
\sum_{n \in E \atop n|w_i} {1 \over n} \geq {1 \over 2} 
\sum_{n \in E \atop n | w_1 \text{ or }w_2} {1 \over n}
& > {1 \over 2} \left ( \sum_{n \in E} {1 \over n} - {2 N^{1-\theta} 
\over (\log\log N)^2 N} \right )\\
& > 1 - O\left ( 
{1 \over N^\theta (\log\log N)^2}\right ).     \endalign
$$
Without loss of generality, assume that the inequality holds for $i=1$, and let 
$E^*$ be those elements of $E$ which divide $w_1$.
\vglue4pt
5.  Use Proposition 2 again, but this time with $J = E^*$,
 $\alpha = \sum_{n \in E^*} 1/n$, and
$\nu = 2/3$, to produce a set $D_j$ satisfying (4.1) and (4.2) with
$E = D_j$.   From (4.5) we have that 
$\Sigma(D_j) < \Sigma(E^*) < (1-e^{-1} + o(1)) \log\log N$.    
\vglue4pt
6.  Set $C_{j+1} = C_j \setminus D_j$.  
If $\sum_{n\in C_{j+1}} {1 \over n} \leq 8/3$, then STOP; 
else, increment $j$ by 1 and go 
back to step 2. 
\medbreak
When this procedure terminates, we are either left with a set 
$D$ from step 3 which proves our proposition, or we are left with
six disjoint sets, $D_1,\ldots ,D_6 \subset
\Cal{C}'(N, N^{1+\delta}; \theta)$ satisfying
$\sum_{n\in D_i} 1/n \in [ 2/3 - 1/ N, 2/3 )$ and 
$$
(e^{-1} - o(1))\log\log N < \Sigma(D_i) < (1 - e^{-1} + o(1))
\log\log N). \tag{4.6}
$$
The lower bound for $\Sigma(D_i)$ follows from Lemma 3 
with $H = D_i$, and the upper bound is
as given in step 5.  

We claim that there exist three of our sets, $D_a, 
D_b, D_c$ such that if $L = 
\Cal{Q}_{D_a} \cap \Cal{Q}_{D_b} \cap \Cal{Q}_{D_c}$, then 
$\Sigma(L) \gg \log\log N$.  For any such triple, we will show that
letting $D = D_a \cup D_b \cup D_c$ satisfies the 
conclusions of Proposition 1.  

To show that $D_a, D_b, D_c$ exist, let $R$ be the 
set of prime powers $\leq N^\theta$ which are contained in at least three 
of the sets $\Cal{Q}_{D_1},\ldots ,\Cal{Q}_{D_6}$.  Then, by (4.6),
$$ \align
\Sigma(R) &> {1 \over 4} \left (
\sum_{ i=1}^6 \Sigma(D_i) \ -\ 
2\sum_{p^a \leq N^\theta \atop p\text{ prime}} {1 \over p^a}
\right ) \\
&> {1 \over 4}\left ( {6 \over e} - 2 - o(1) \right ) \log\log N
\gg \log\log N. \endalign
$$
Thus, since there $20 = {6\choose 3}$ triples of sets chosen from
$\{D_1,\ldots ,D_6\}$, there is at least one such triple which
gives $\Sigma(L) > \Sigma(R)/20 \gg \log\log N$. 

Now, letting $D = D_a \cup D_b \cup D_c$ 
certainly satisfies (1.2).  Suppose that the number of elements of  
$D$ which do not divide any element of $I$ is at most 
$N^{1-\theta}/(\log\log N)^2$. 
Then, the hypotheses of Proposition~3 hold for $E = D_a, D_b$,
and $D_c$.  Case B of Proposition 3 cannot hold for  
$E = D_a$ (or $D_b,$ or $D_c$), else (4.5) and (4.6) would give us 
$$
\eqalign{
\sum_{q | \text{\rm gcd}(w_1,w_2)\atop q \in \Cal{Q}_{D_a}} 
{1 \over q} &> \sum_{q | w_1\atop q\in \Cal{Q}_{D_a}} {1 \over q}
+ \sum_{q|w_2\atop q\in \Cal{Q}_{D_a}} {1 \over q} - 
\sum_{q \in \Cal{Q}_{D_a}} {1 \over q} \cr
&> \left ( {3 \over e} - 1 - o(1) \right ) \log\log N \gg \log\log N,
}
$$
which, by Lemma 2, would imply that $w_1 = w_2$.
Thus, case A of Proposition~3 holds for $E = D_a, D_b$, and $D_c$:  
Let $W_a$, $W_b$, and $W_c$ 
be the single integer in $I$ dividing all elements of $D_a, D_b,$ and $D_c$, 
respectively, and thus they are all divisible by every element of $L$. 
Since $\Sigma(L) \gg \log\log N$, we have, from Lemma~2, that 
$W_a = W_b = W_c = W$, for some $W \in I$.  Proposition 1 now follows
since $\text{lcm} \{ n \in D \} | W$.

\section{Proof of Proposition 2}

To prove Proposition 2 we will need the following lemma.

\nonumproclaim{Lemma 4} Suppose $S$ is a set of integers{\rm ,} all of whose prime
power divisors are less than $N${\rm ,} which satisfies 
$\sum_{n \in S} 1/n \geq \rho > \mu$.  If $N$ is large in terms of
$\rho$ and $\mu${\rm ,} then there exists a subset $T \subseteq S$ for which 
$$
\sum_{n \in T} {1 \over n} > \mu, \text{ and } 
\sum_{n \in T \atop q | n} {1 \over n} > {\rho - \mu \over 2q
\log\log{N}},
\ \text{ for all $q \in \Cal{Q}_T$}. \tag{5.1} 
$$
\endproclaim

\demo{Proof}  We form a chain of subsets $S_0 := S 
\supset S_1 \supset \cdots \supset T := S_k$ as 
follows:  given $S_i$, let $q_i$ be the smallest prime power 
such that 
$$
\sum_{n \in S_i\atop q_i|n} {1 \over n} < 
{\rho - \mu \over 2q_i\log\log N},
$$
if such $q_i$ exists, and then let $S_{i+1} = S_i
\setminus \{ n \in S_i : q_i | n\}$.  If no such $q_i$ exists,
then let $k = i$ and $T = S_i = S_k$.
We have that
$$
\sum_{n \in T} {1 \over n} > \rho - {\rho - \mu \over 
2\log\log{N}} \sum_{p^a \leq N \atop p \text{ prime}} { 1 \over p^a}  
> \mu,
$$
for $N$ large enough, since $\sum_{p^a \leq N} 1/p^a < 2\log\log N$. 
\enddemo

\demo{Proof of Proposition {\rm 2}} 
We first use Lemma 4 with $\rho = \alpha$, $\mu = \nu$, 
and $S = J$, to produce a set 
$D_0 = T$ satisfying (5.1). 
Thus, (4.2) holds for $E = D_0$.

We will construct a chain of subsets 
$D_0\supset D_1\supset D_2\supset \cdots$,
where each set $D_j$ satisfies (4.2) with $E = D_j =
D_{j-1} \setminus \{w_j\}$, where $w_j$ is
some yet to be chosen element of $D_{j-1}$.
If we can do this then we will eventually reach a set $D_k$
which also satisfies (4.1), since each $w_j \geq N$, and so the proposition 
will be proved.

Suppose (4.2) is satisfied for $E = D_{j-1}$, for $j \geq 1$.  
Take Lemma 4 with $S = D_{j-1}$, $\rho = \nu$, and $\mu = \nu/2$,
and let $w_j$ be the smallest element of $T$.  Let $q \in \Cal{Q}_
{D_j}$.  If $q \nmid w_j$, then
$$
\sum_{n \in D_j \atop q | n} {1 \over n} = \sum_{n \in D_{j-1}
\atop q | n} {1 \over n} > {\min\{\nu, \alpha -\nu\} \over 5q \log\log N},
$$
by hypothesis.  On the other hand, if $q | w_j$, then, by (5.1), we get
$$
\sum_{n \in D_j \atop q|n} {1 \over n} \geq 
\sum_{n \in T \atop q|n} {1 \over n} - {1 \over w_j}
> {\nu \over 4q\log\log N}\ -\ {1 \over N} > {\nu
\over 5q \log\log N},
$$
since $q \leq N^\theta$, with $\theta < 1$, and $\nu \gg 1$, and so  
(4.2) holds for $E = D_j$.
\enddemo

\section{Proof of Proposition 3}

Let $E_I$ denote the set of integers 
in $E$ which divide an integer in $I$.  Then we have, by hypothesis, 
that $|E_I| > |E| - N^{1-\theta}/(\log\log N)^2$.  
If $q \in \Cal{Q}_E$, then
$$
\sum_{n \in E_I \atop q | n} {1 \over n} > \sum_{n \in E
\atop q|n} {1 \over n}\ -\ {N^{1-\theta} \over N (\log\log N)^2} \gg
{1 \over q \log\log N}, \tag{6.1}
$$
since $q \leq N^\theta$ and $E$ satisfies (4.2).  Thus, we have
that $\Cal{Q}_{E_I} = \Cal{Q}_E$. 

We will show at the end of this section that
for all $q \in \Cal{Q}_E$, there exists an 
integer $qd \in [N^{3/4} , N^{3/4+\theta}]$ such that 
$$
\sum_{n \in E_I \atop qd|n}\ {1 \over n} \gg_\theta {1 \over qd 
(\log\log N)^2},\tag{6.2}
$$
where $\omega(d) \leq \omega_0 = \log\log N/\log\log\log\log N$, for
$N$ sufficiently large, and all the 
prime divisors of $d$ are greater than 
$y := \exp((1/8 - \theta/2)\log N/\log\log N)$.

For now, let us assume that this is true and let $qd$ satisfy (6.2) for 
a given $q \in \Cal{Q}_E$.
All the elements of $E_I$ which are divisible by $qd$ must divide
the same number $n(q) \in I$, since otherwise there are two distinct 
numbers $n_1(q)$ and $n_2(q)$ which differ by $\leq N^{3/4}$ but yet
are both divisible by $qd > N^{3/4}$, which is impossible.
We will show that as a consequence of this and (6.2), 
$$
\sum_{p^a | n(q)\atop p^a \in \Cal{Q}_E} {1 \over p^a} > 
\left ({1 \over e} - o(1) \right )\log\log N. \tag{6.3}
$$
This implies there are at most two distinct values of $n(q)$,
for all $q \in \Cal{Q}_E$:  for if there were three prime powers $q_1, 
q_2, q_3$ with $n(q_1), n(q_2), n(q_3)$ distinct, then, by Lemma 2,
$$
\sum_{p^a | \text{gcd}(w_1,w_2)} {1 \over p^a} \ll \log\log\log N,
$$
so that, by (6.3),  
$$
\align
\log\log N + O(1) = \sum_{p^a \leq N \atop \text{$p$ prime}} {1 \over p^a} 
> \Sigma(E) &\geq \sum_{i=1}^3 \sum_{p^a | n(q_i) \atop
p^a \in \Cal{Q}_E} {1 \over p^a} + O(\log\log\log N)\\
&> (3e^{-1} -o(1))\log\log N,
\endalign 
$$
which is impossible.

If there is just one value for $n(q)$, for all $q \in \Cal{Q}_
E$, then $w = n(q)$ satisfies case A of Proposition 3:
Otherwise, there are two possible values for $n(q)$, call them $w_1$
and $w_2$, which satisfy (4.4).  The lower bound in (4.5)
comes from (6.3).  Moreover,   
$$
\sum_{q | w_1\atop q \in \Cal{Q}_E} {1 \over q} \leq \sum_{p^a \leq N}
{1 \over p^a} - \sum_{q | w_2 \atop q \in \Cal{Q}_E} {1 \over q}
+ \sum_{p^a | \text{gcd}(w_1,w_2)} {1 \over p^a} \leq
(1 - e^{-1} + o(1))\log\log N,
$$
which implies the upper bound in (4.5) (note: the same upper bound holds 
for $w_2$), using the Prime Number Theorem (6.3), and Lemma 2, respectively. 

If $w_1, w_2$ fail to satisfy (4.3), then 
$$ \align
\#\{  n \in E_I:n \nmid w_1\text{ or }
w_2\}& > \#\{ n \in E : n \nmid w_1\text{ or }w_2\} - {N^{1-\theta}
\over (\log\log N)^2}\\
& > {N^{1-\theta} \over (\log\log N)^2}. \endalign
$$
Since there are $\leq N^{3/4}$ integers in $I$, there must exist an integer 
$x \in I, x \neq w_1$ or $w_2$, for which
$$
\#\{ n \in E_I :  n | x\} \gg {N^{1-\theta} \over N^{3/4}(\log\log N)^2}
 = {N^{1/4-\theta}\over (\log\log N)^2}. \tag{6.4}
$$
Therefore,
$$ \align
\text{lcm}_{n \in E, n|x}\  n &\leq \text{gcd}(x,w_1w_2)\\& \leq 
\text{gcd}(x,w_1)\text{gcd}(x,w_2) < (x-w_1)(x-w_2) < N^{3/2}; \endalign
$$
but then we have
$$
\#\{n \in E\ :\ n | x\} \leq \tau \left ( \text{lcm}_{n \in E,
n | x}\ n\ \right ) \leq \max_{l \leq N^{3/2}} \tau(l) = N^{o(1)},
$$
which contradicts (6.4), and so (4.3) follows.  Thus, the proof of Proposition
3 is complete once we establish (6.2) and (6.3).

To show (6.3), we observe that every integer $m \in F = 
\{ n/qd : n \in E$, $qd | n\}$ satisfies $\omega(m) \sim \log\log N$, 
since $\omega(qd) \leq \omega_0 = o(\log\log N)$, and since $E \subset
\Cal{C}'(N, N^{1+\delta}; \theta)$.  From this and (6.2), 
$F$ satisfies the hypotheses of Lemma~3 with $H = F$.  
Thus, $\Sigma(F) > (e^{-1} - o(1))\log\log N$, which implies (6.3).

We will now establish (6.2).  First, we claim that for every 
$n \in E$, where $q|n$ and $q\in \Cal{Q}_E$, there exists a 
divisor $qd \in [N^{3/4},N^{3/4+\theta}]$,
where $p | d$ implies $p > y$ (though it may not be the case that
$\omega(d) \leq \omega_0$).  To show this, we construct such a $d$
by adding on prime factors one at a time, until $qd$ is in this interval.
There are enough prime factors $> y$ to do this, since for $N \gg_\epsilon 1$
we have
$$
\prod_{p^a || n/q \atop p > y} p^a > {n \over q \prod_{p^a || n \atop p\leq y}
p^a} > {n \over q y^{\Omega(n)}} > {N \over N^\theta \exp \left (
{(1/8-\theta/2)\Omega(n) \log N \over \log\log N} \right )} >  N^{3/4},
$$
for $N$ sufficiently large, since $\Omega(n) \sim \log\log N$.

If (6.2) fails to hold for all $d \in [N^{3/4}/q, N^{3/4+\theta}/q]$ with
$\omega(d) \leq \omega_0$, then we would have by (4.2) and Mertens' theorem 
that, 
$$
\align 
\noalign{\vskip-6pt}
&\tag{6.5}\\
{\min\{\nu, \alpha - \nu\} \over 5q \log\log N} &< 
\sum_{n \in E\atop q|n} {1 \over n} < \sum_{N^{3/4}/q \leq d \leq 
N^{3/4+\theta}/q \atop p|d \Rightarrow p > y} \sum_{n \in E \atop qd|n} 
{1 \over n} 
\\
&< \sum \Sb N^{3/4}/q \leq d \leq N^{3/4+\theta}/q \\
            p|d \Rightarrow p > y \\
            \omega(d) < \omega_0 \endSb 
\sum_{n \in E \atop qd|n} {1 \over n} +
   \sum_{d: p|d \Rightarrow y < p < N\atop \omega(d) \geq \omega_0} 
\sum_{m \leq N^{1+\delta}/qd\atop (n=mqd)} {1 \over qdm}\cr
& = o\left ( {1 \over q(\log\log N)^2} \sum_{d :  p|d \Rightarrow
y < p < N} {1 \over d} \right )\! \!+\! O\left ( {\log N \over q}\! \sum_{d : p|d 
\Rightarrow y < p < N \atop w(d) \geq \omega_0} {1 \over d} \right )\!.\\
\endalign
$$ 
Now,
$$
\sum_{d :  p|d \Rightarrow y < p < N} {1 \over d} \leq 
\prod_{y < p < N\atop \text{$p$ prime}} \left (1 - {1 \over p} \right)^{-1}
 \ll {\log N \over \log y} \ll \log\log N,
$$
by Mertens' theorem, and for $k = (\log\log\log N)^3$, we have, again by 
Mertens' theorem,  
$$
\eqalign{
\sum_{d : p|d \Rightarrow y < p < N \atop \omega(d) \geq \omega_0} {1 \over d}
&\ll \sum_{d : p|d \Rightarrow y < p< N} {k^{\omega(d) - \omega_0} \over d}
= {1 \over k^{\omega_0}} \prod_{y < p < N\atop p\text{ prime}} \left ( 1 +
{k \over p - 1} \right )\cr
& = {1 \over k^{\omega_0}} \left ( {\log N \over 
\log y} \right )^{k + o(k)} \ll {1 \over \log^2 N}.
}
$$ 
Combining these two applications of Mertens' theorem with (6.5), we arrive
at a contradiction.  Thus, there must exist a $d \in 
[N^{3/4}/q,N^{3/4+\theta}/q]$ satisfying (6.2), with 
$\omega(d) \leq \omega_0 = o(\log\log N)$. 

\section{Acknowledgements}

First, and foremost, I would like to thank my advisor Andrew Granville for
his encouragement and for helping me edit this paper to get it into its 
current form.  I would also like to thank P. Erd\H os and R. L. Graham for the 
wonderful questions.

\references

[1]
\name{K.\ Dickman}, On the frequency of numbers containing prime factors
of a certain relative magnitude, {\it Arkiv.\ Math.\ Astr.\ Fys\/}.\
{\bf 22} (1930), 1--14.

[2]
 \name{P. Erd\H os} and \name{R.\ L.\ Graham}, Old and new problems
and results in combinatorial number theory, {\it Enseign.\ Math\/}.\
(1980), 30--44.

[3]
\name{R.\ K.\ Guy}, {\it Unsolved Problems in Number Theory\/}, Second
edition, Springer-Verlag, New York, 1994, 158--166.

[4
\name{H.\ Halberstam} and \name{H.-E.\ Richert}, {\it Sieve Methods\/},
{\it London Math.\ Soc.\ Monographs\/}, No.\ 4  (1974), Academic Press,
New York, 1974.

[5]
\name{H.\ L.\ Montgomery}, Ten lectures on the interface
between analytic number theory and harmonic analysis, {\it C.B.M.S.\
Reg.\ Conference Series Math\/}.\ {\bf 84}, Amer.\ Math.\ Soc.,
Providence, RI, 1994.

\endreferences
\bye

\centerline{(Received May 16, 2001)}
\enddocument